\documentclass[12pt]{article}
\usepackage{amsfonts,amsmath,bbm}
\def\proof{\noindent{\bf Proof:}\hskip10pt}        
\def\QED{\hfill $\Box$}

\font\tenmath=msbm10 scaled 1200
\font\sevenmath=msbm7 scaled 1200
\font\fivemath=msbm5 scaled 1200
\newfam\mathfam \textfont\mathfam=\tenmath
\scriptfont\mathfam=\sevenmath \scriptscriptfont\mathfam=\fivemath

\begin{document}
\def \\ { \cr }
\def \R{\mathbb{R}}
\def\N{\mathbb{N}}
\def\E{\mathbb{E}}
\def\P{\mathbb{P}}
\def\Z{\mathbb{Z}}
\def\D{\mathbb{D}}
\def\C{\mathbb{C}}
\def\da{^{\downarrow}}
\def \e{{\rm e}}
\def \p{{\cal P}}
\def \s{{\cal S}}
\def \g{{\cal G}}
\newtheorem{theorem}{Theorem}
\newtheorem{definition}{Definition}
\newtheorem{proposition}{Proposition}
\newtheorem{lemma}{Lemma}
\newtheorem{corollary}{Corollary}
\centerline{\LARGE \bf Reflecting a Langevin Process}
\vskip 2mm
\centerline{\LARGE \bf    at an Absorbing Boundary}

\vskip 1cm
\centerline{\Large \bf Jean Bertoin}
\vskip 1cm
\noindent
\centerline{\sl Laboratoire de Probabilit\'es et Mod\`eles Al\'eatoires}
\centerline{\sl Universit\'e Pierre et Marie Curie (Paris 6), et C.N.R.S. UMR 7599}
\centerline{\sl 175, rue du Chevaleret} 
\centerline{\sl F-75013 Paris, France}
\vskip 15mm

\noindent{\bf Summary. }{\small  We consider a Langevin process with white noise random forcing. We suppose that the energy  of the particle is instantaneously absorbed when it hits some fixed obstacle. We show that nonetheless, the particle can 
be instantaneously reflected, and study some properties of this reflecting solution.}
 \vskip 5mm
\noindent
{\bf A.M.S. Classification.}  {\tt Primary 60 J 50, 60 J 55. Secondary 60 G 18, 60 G 52}
\vskip 3mm
\noindent{\bf e-mail.} {\tt jbe@ccr.jussieu.fr}

\begin{section}{Introduction}
Langevin \cite{Langevin} introduced a probabilistic model to describe the evolution of a particle
under random external forcing.
Langevin's motion has smooth trajectories, in the sense that if 
$Y_t$ stands for the position of the particle at time $t$, then
the velocity $\dot Y_t:={\rm d}Y_t/{\rm d}t$ 
is  well-defined and finite everywhere.
We assume that the external force (i.e. the derivative of the velocity)  is a white noise, so that
$$ \dot Y_t=  \dot Y_0 + W_tÊ\,,$$
 where $W=(W_t, t\geq 0)$ is a standard Wiener process\footnote{
 We mention that Langevin considered more generally the case when 
 the velocity is given by an Ornstein-Uhlenbeck process.}, and $Y_t$ can then be expressed in terms of the integral of the latter. Plainly the Langevin process $Y=(Y_t,t\geq 0)$ is not Markovian, however
 the pair $(Y,\dot Y)$, which is often called Kolmogorov's process, is Markov.
 We refer to Lachal \cite{LaSP} for an interesting  introduction to this field,  historical comments,  and a long list of references.
 
Recently,  Maury \cite{Maury}  has considered the situation when the particle may hit an obstacle
and   then its velocity is instantaneously absorbed.  Roughly, this arises in a model  for 
the motion of individuals during the formation of a crowd.
For the sake of simplicity,
we shall work in dimension $1$ and assume that the particle evolves in the half-line
$[0,\infty[$. As long as the particle remains in $]0,\infty[$, its dynamics are those described above. 
When the particle reaches the boundary point $0$, we suppose that its energy  is 
absorbed, in the sense that the velocity instantaneously drops to $0$.

There is an obvious way to construct a process that fulfills these requirements. Namely,
given an initial position $x\geq 0$ and an initial velocity $v\in\R$, if we consider the free Langevin process  (i.e. when there is no obstacle) started from $x$ with velocity $v$ 
$$Y_t= x+\int_{0}^{t} (v+W_s){\rm d}s $$ 
and define
$$\zeta=\inf\{t\geq 0: Y_t=0\}\,,$$
then we can take $(Y_{t\wedge \zeta}, t\geq 0)$.
In the sequel, we shall refer to this solution as the Langevin process stopped
at $0$.

It is easy to show that when $x=v=0$, the free Langevin process $Y=
(Y_t, t\geq 0)$ returns to $0$ at arbitrarily small times, and that its velocity at such return times is never $0$. In the presence of an obstacle, this suggests that reseting  the velocity to $0$ at every hitting time of the boundary might impede the particle to ever exit $0$, so that the Langevin process stopped at $0$ might be the unique solution.

The purpose of this work is to point at the rather surprising fact (at least for the author) that there  exists a remarkable reflecting  solution which exits instantaneously  from the boundary point $0$.
Constructing recurrent extensions
of a given Markov process killed when it reaches some boundary 
is a classical problem in the theory of Markov processes;
see in particular It\^o \cite{Ito}, Rogers \cite{Rogers1, Rogers2}, Salisbury \cite{Salisbury}, Blumenthal \cite{Bl} and references therein
(note that in these works, the boundary is reduced to a single point
whereas here, the Markov process $(Y,\dot Y)$ is killed when it reaches
$\{0\}\times \R$). However it does not seem easy to apply directly this theory
 in our setting, as technically, it would require the construction of a so-called 
 entrance law from $(0,0)$  for the killed Kolmogorov process.
 Some explicit information about the transitions of latter
 are available in the literature, see in particular Lachal  \cite{La91, LaSP} and references therein. Unfortunately, the expressions are quite intricate and tedious to manipulate.
 
 We shall first construct a reflecting solution by combining the classical idea of Skorohod 
 with techniques of time substitution. Then,  we shall establish
 uniqueness in distribution. In this direction, we have to obtain {\it a priori} information about the measure of the excursions away from $0$ for any reflected solution. More precisely, we will first
 observe that the velocity of a reflecting Langevin process does not remain strictly positive immediately  after the beginning of an excursion away from $0$,
 which may be a rather surprising fact. Then, combining a fundamental  connection between the Langevin process and the symmetric stable process of index $1/3$ with a result of Rogozin \cite{Rogozin} on the two-sided exit problem 
for stable L\'evy processes,  we shall compute explicitly certain distributions
under the excursion measure of any reflecting solution. This will enable us to characterize the excursion measure and hence to establish uniqueness in law.
 Finally, we turn our attention to statistical self-similarity properties of the reflected Langevin process; in particular we shall
 establish that the Hausdorff dimension of the set of times at which the particle reaches the boundary is $1/4$ a.s.

\end{section}

\begin{section}{Construction of the reflecting process}
The purpose of this section is to provide an explicit construction of a Langevin process with white noise external forcing, such that the energy of the particle is absorbed each time it hits  the boundary point $0$,  and which is nonetheless instantaneously reflected at $0$. In order to focus on the most interesting situation, we shall suppose thereafter that the initial location and velocity are both zero, although the construction obviously extends to arbitrary initial conditions.

By this, we mean that we are looking for a c\`adl\`ag process 
$(X,V)=((X_t, V_t), t\geq 0)$
with values in $[0,\infty[\times \R$, which starts from $X_0=V_0=0$ 
and fulfills the following five requirements. First, 
 $V$ is the velocity process of $X$, that is
\begin{equation}\label{Cond1}
X_t\,=\, \int_0^t V_s {\rm d}s\,.
\end{equation}
Second, we require that the energy of the particle is absorbed at the boundary, viz.
\begin{equation}\label{Cond2}
V_t=0\qquad \hbox{for every $t\geq 0$ such that }X_t=0, \quad \hbox{a.s.},
\end{equation}
and, third,  that the process $X$ spends zero time at $0$ :
\begin{equation}\label{Cond3}
 \int_0^{\infty} {\bf 1}_{\{X_t=0\}}{\rm d}t\,=\, 0\qquad \hbox{a.s.}
 \end{equation}
Fourth, the process $X$ evolves  in $]0,\infty[$ as a free Langevin process, that is
 if we define for every $t\geq 0$ the first return time to the boundary after $t$, 
 $$\zeta_t=\inf\{s\geq t: X_{s}=0\}\,,$$
  then for every stopping time $S$ in the natural filtration
  (after the usual completions)  $({\mathcal F}_t)_{t\geq 0}$ of $X$ :
\begin{eqnarray}\label{Cond4}
& &
\hbox{conditionally on $X_S=x>0$ and $V_S=v$,  the process
$(X_{(S+t)\wedge \zeta_S}, t\geq 0)$}\nonumber \\
& &\hbox{is independent of $ {\mathcal F}_S$, and has the distribution of  
a Langevin process }\nonumber \\
& &\hbox{started with velocity 
$v$  from the location $x$ and stopped at $0$.}
\end{eqnarray}
 Roughly, the final condition is a natural requirement of  regeneration at return times to the boundary. Specifically,  for every stopping time $S$ such that $X_S=0$ a.s., we have that
 \begin{equation}\label{Cond5}
\hbox{the process  $(X_{S+t}, t\geq 0)$ has the same law as $(X_t, t\geq 0)$ and is independent of ${\mathcal F}_{S}$.}
\end{equation}
Conditions (\ref{Cond4}) and (\ref{Cond5}) entail that $(X,V)$ enjoys the strong Markov property; note however the latter has jumps at predictable stopping times and thus is not a Hunt process.
Note also that, since a free Langevin process started from an arbitrary position with an arbitrary velocity eventually reaches the boundary point $0$ a.s., $(0,0)$ is necessarily a recurrent point for $(X,V)$.

Let $W=(W_t,t\geq 0)$ be a standard Wiener process started from $W_0=0$. 
Recall that $W$ possesses a local time  at $0$ which is defined by
$$L_t:=\lim_{\varepsilon\to 0+}\frac{1}{\varepsilon}\int_0^t
{\bf 1}_{\{0<W_s<\varepsilon\}}{\rm d}s\,,Ê\qquad t\geq 0\,.$$
The process  $L=(L_t, t\geq 0)$ has continuous paths and the support
of the Stieltjes measure ${\rm d}L_t$ coincides with
the set of times $t$ at which $W_t=0$, a.s.
We write
$$\tau_{t}:=\inf\{s\geq 0: L_s>t\}\,, \qquad t\geq 0$$
for the right-continuous inverse of $L$.

Next, we define the free  
Langevin process
$$Y_t\,=\,\int_0^{t} W_s{\rm d}s\,,\qquad t\geq 0\,.$$
The scaling and strong Markov properties of the Wiener process
readily entail the following lemma which will be crucial for our analysis.

\begin{lemma}\label{Lemma1} The time-changed process
$$\sigma_t:=Y_{\tau_t}\,, \qquad t\geq 0$$
is a symmetric stable L\'evy process with index $1/3$.
Moreover, $Y_s$ lies between $\sigma_t$ and $\sigma_{t-}$ for every 
$s\in[\tau_{t-},\tau_t]$.
\end{lemma}

 We refer to Biane and Yor \cite{BY} for the first claim. The second
follows from the fact that $Y$ is monotone on every time-interval $[\tau_{t-},\tau_t]$
on which the Wiener process $W$ makes an excursion away from $0$. 
It will play an important part in our study as it enables us to relate exit properties
of the Langevin process to that of the symmetric stable process $\sigma$, and useful information on the latter can be gleaned from the literature.

We also consider the  non-increasing infimum process
$$I_t=\inf\{Y_s : 0\leq s \leq t\}\,.$$
Plainly,  $I$ is absolutely continuous (because so is $Y$)
and ${\rm d}I_t=0$ whenever $Y_t>I_t$.
We observe that $W_t\leq 0$ for every $t$ such that $Y_t=I_t$. Indeed, if we had
$W_t>0$ for such a time $t$, then $Y$ would be strictly increasing on some neighborhood of
$t$, which is incompatible with the requirement that $Y_t=I_t$. 
Now for every $t\geq 0$ such that $Y_t=I_t$ and $W_t<0$, $Y$ is strictly decreasing
on some interval $[t,t']$ with $t'>t$, and thus $Y=I$ on $[t,t']$. Since the total time that $W$ spends at $0$ is zero, we conclude that
$$I_t\,=\,\int_0^t {\bf 1}_{\{Y_s=I_s\}}W_s {\rm d}s\,.$$

Next, consider the process $\tilde X_t=Y_t-I_t$ which can thus be expressed in the form
$$\tilde X_t=\int_0^t \tilde V_s {\rm d}s\,,$$
where
$$\tilde V_t={\bf 1}_{\{Y_t>I_t\}}W_t= {\bf 1}_{\{\tilde X_t>0\}}W_t \,.$$
We observe that the velocity process $\tilde V=(\tilde V_t, t\geq 0)$ has c\`adl\`ag paths a.s.
Indeed, if $\tilde V_t\neq 0$, then $\tilde X_t>0$, so the indicator function
$s\to {\bf 1}_{\{\tilde X_s>0\}}$ remains equal to $1$ on some neighborhood of $t$ and 
$s\to {\bf 1}_{\{\tilde X_s>0\}} W_s$ is thus continuous at $t$. Plainly, $s\to {\bf 1}_{\{\tilde X_s>0\}}W_s$ is also continuous at $t$ when $W_t=0$. Suppose now that $\tilde V_t= 0$ and $W_t\neq 0$, and recall that necessarily $W_t<0$. Then,  $Y$ is strictly decreasing on some neighborhood of $t$, say $]t-\varepsilon, t+\varepsilon[$. As $Y_t=I_t$, either
$ {\bf 1}_{\{\tilde X_s>0\}}= 1$ for $t-\varepsilon<s<t$ and $ {\bf 1}_{\{\tilde X_s>0\}}= 0$ for $t\leq s<t+\varepsilon$, or    the indicator function $ {\bf 1}_{\{\tilde X_s>0\}}$ remains equal to $0$
 on some neighborhood of $t$. In both cases, this implies that $s\to {\bf 1}_{\{\tilde X_s>0\}} W_s$ is c\`adl\`ag  at $t$.

In conclusion, $(\tilde X, \tilde V)$ is a c\`adl\`ag process which fulfills the requirements
(\ref{Cond1}) and (\ref{Cond2}), and an immediate application of the Markov property
shows that (\ref{Cond4}) also holds. However, (\ref{Cond3}) clearly fails, as the closed set
$${\mathcal I}:=\{t\geq0: \tilde X_t=0\}$$
has a strictly positive (as a matter of fact, infinite) Lebesgue measure.
This incites us to study ${\mathcal I}$ in further details, and in this direction, the following notion
appears naturally. For every $t\geq 0$, we define the random set  
$${\mathcal J}:=\left\{t\geq 0: W_t<0, \tilde X_t=0 \hbox{ and } \tilde X_{t-\varepsilon}
>0\hbox{ for all $\epsilon>0$ sufficiently small}\right\}.$$
In other words, ${\mathcal J}$ is the set of times at which $Y$ reaches its infimum
for the first time during a negative excursion of $W$.
Observe that each point $s\in{\mathcal J}$ is necessarily isolated in ${\mathcal J}$, and in particular  ${\mathcal J} $ is countable.
Further, we introduce the notation
$$d_t:=\inf\{s>t: W_s=0\}$$
for the first return time of the Wiener process $W$  to $0$ after time $t$.
We may now
state the following lemma.

\begin{lemma}\label{Lemma2} 
The canonical decomposition of the interior ${\mathcal I}^{o}$ of ${\mathcal I}$ as the
union of disjoint open intervals is given by
$${\mathcal I}^{o}\,=\,\bigcup_{s\in{\mathcal J}}]s,d_s[\,.$$
Moreover, the boundary $\partial {\mathcal I}={\mathcal I}\backslash 
{\mathcal I}^{o}$ has zero Lebesgue measure. 
\end{lemma}

\proof We start by recalling  that the free Langevin process $Y$ oscillates at the initial time,
in the sense that
$$\inf\{t>0: Y_t>0\}=\inf\{t>0: Y_t<0\}=0\qquad \hbox{a.s.}$$
Indeed, for every $t>0$, $Y_t$ is a centered Gaussian variable,
so $\P(Y_t>0)=\P(Y_t<0)=1/2$, and the claim follows from Blumenthal's $0$-$1$ law for the Wiener process. 

Now pick an arbitrary rational number $a>0$ and work conditionally on
$a\in{\mathcal I}^{o}$. Then we know that $W_a\leq 0$, and as
 $\P(W_a=0)=0$,  we may thus assume that
 $W_a< 0$. Recall that $d_a=\inf\{s>a: W_s=0\}$ denotes the first hitting time of $0$ by $W$ after time $a$. Plainly $Y$ decreases on $[a,d_a]$, and thus $]a,d_a[\subseteq {\mathcal I}^{o}$.
  An application of the strong Markov property of the Wiener process $W$ at the stopping time $d_a$, combined with the oscillation property of the free Langevin process mentioned above, shows that
 $d_a$ is the right extremity of the interval component of  ${\mathcal I}^{o}$ which contains $a$.

Next, consider $g_a=\sup\{t<a:W_t=0\}$, the left-extremity of the interval containing $a$ on which
$W$ makes an excursion away from $0$. We shall show that $]g_a,d_a[$ intersects ${\mathcal J}$
at a single point $s$. Recall that $L=(L_t, t\geq 0)$ denotes the process of the local time of $W$ at $0$
and $\tau_t=\inf\{s\geq 0: L_s>t\}$ its right-continuous inverse, and that the compound process
$\sigma:=Y\circ \tau$ is a symmetric stable L\'evy process with index $1/3$
(see Lemma \ref{Lemma1}). We write $\iota_t:=\inf_{0\leq s\leq t}\sigma_s$ for the infimum process of $\sigma$.
The excursion interval $]g_a,d_a[$ corresponds to a jump of the inverse local time,
that is $]g_a,d_a[=]\tau_{t-},\tau_t[$ for some $t>0$.
 Because $Y$ is monotone on every excursion interval of $W$,
the identities
$$\iota_v\,=\,I(\tau_v)\quad \hbox{and}\quad
\iota_{v-}\,=\,I( \tau_{v-})$$
hold for every $v$, a.s.
Since $Y_{d_a}=I_{d_a}$, the stable process $\sigma$ reaches a new infimum at time $t$.
On the other hand, we know that a symmetric stable process never jumps at times $v$
such that $\sigma_{v-}=\iota_{v-}$,
see for instance Rogers \cite{Rogers}, so
$$I_{ \tau_{t-}}=\iota_{t-}\,<\sigma_{t-}=Y_{\tau_{t-}}\,.$$ 
That is  $I_{g_a}<Y_{g_a}$, and
as the Langevin process $Y$ decreases on $]g_a,d_a[$ and $Y_a=I_a$, this
implies that  $s:=\inf\{u>g_a: Y_u=I_u\}\in{\mathcal J}$. More precisely, we have 
that $]g_a,d_a[\cap {\mathcal J}=\{s\}$, and it should now be plain that
the interval component of ${\mathcal I}^{o}$ which contains $a$ is $]s,d_a[=]s,d_s[$ a.s.
Taking the union over the set of rational $a$'s with $a\in{\mathcal I}^o$ establishes the inclusion
$${\mathcal I}^o\,\subseteq \bigcup_{s\in{\mathcal J}}]s,d_s[\qquad \hbox{a.s.},$$
and the converse inclusion is obvious.

Finally, we turn our attention to the second assertion. Pick any $s\in \partial {\mathcal I}$; we then know that $W_s\leq 0$. Note that if $W_s<0$, then  $s\in{\mathcal J}$.
Since the set of times at which the Wiener process is $0$ has zero Lebesgue measure
and ${\mathcal J}$ is countable, this shows the second claim. \QED

We are now able to complete the construction. 
Introduce 
$$T_t:=\inf\left\{s\geq 0: \int_0^s {\bf 1}_{\{\tilde X_u>0\}}{\rm d}u>t\right\}\,, \qquad t\geq 0\,,$$
and set
$$X_t=\tilde X\circ T_t \qquad \hbox{and}\qquad  V_t=\tilde V\circ T_t \qquad \hbox{for every }t\geq 0\,.$$
 Informally, the time-substitution
amounts to suppressing the pieces of the paths on which the free Langevin process is at its minimum,
and this makes the following statement rather intuitive.
Observe that the process  $t\to T_t$ is  increasing and right-continuous, so the 
sample paths of $((X_t,V_t), t\geq 0)$ are c\`adl\`ag a.s.

\begin{theorem}\label{Thm1} The process $(X,V)$ constructed above fulfills the requirements
(\ref{Cond1}), (\ref{Cond2}),  (\ref{Cond3}), (\ref{Cond4}) and (\ref{Cond5}).
\end{theorem}

\proof  The argument relies essentially on a change of variables formula
for the increasing process $T$ that we shall first justify with details.

For every $t\geq 0$, we consider the partition of the interval $[0,T_t]$ into three Borel sets
$$B_1(t):=[0,T_t]\cap {\mathcal I}^c\quad ,\quad B_2(t)
:=[0,T_t]\cap \partial {\mathcal I}\quad ,\quad 
B_3(t):=[0,T_t]\cap {\mathcal I}^o\,$$
where ${\mathcal I}^c=\{s\geq 0: \tilde X_s>0\}$.
Recall Lemma \ref{Lemma2}, and in particular that $\partial {\mathcal I}$ has zero Lebesgue measure.
We deduce that the canonical decomposition of the increasing process $T$ into  its continuous component and its jump component is
$$T_t\,=\,\int_0^{T_t}{\bf 1}_{\{\tilde X_s>0\}}{\rm d}s +  \sum_{s\in{\mathcal J}, s\leq T_t}(d_s-s)
\,=\, t+\sum_{s\in{\mathcal J}, s\leq T_t}(d_s-s)\,,$$
where the second equality stems from the very  definition of the time-substitution.

The classical change of variables formula gives
\begin{equation}\label{eqn6}
f(T_t)\,=\, f(0)+ \int_0^t f'(T_s){\rm d}s
+ \sum_{s\in{\mathcal J}, s\leq T_t}(f(d_s)-f(s))\,,
\end{equation}
where $f:\R_+\to \R$ stan{\rm d}s for a generic function which is continuously differentiable.
A standard argument using the monotone class theorem shows that (\ref{eqn6})
extends  to the case
when $f$ has a derivative $f'$ in Lebesgue's sense which is locally bounded.
Apply the change of variables formula for 
$$f(s)\,=\,\tilde X_s\,=\,\int_0^s \tilde V_u {\rm d}u\quad , \quad f'(s)=\tilde V_s\,.$$
Since $f(d_s)-f(s)=0$ for every $s\in{\mathcal J}$, 
this gives
$$X_t=f(T_t)=\, \int_0^t   \tilde V_{T_s}{\rm d}s=\, \int_0^t   V_s{\rm d}s\,,$$
which proves (\ref{Cond1}).

 Next,
as the pair of processes $(\tilde X, \tilde V)$ fulfills (\ref{Cond2}), the time-substitution
 shows that the same holds for $(X,V)$.
Then, we apply the change of variables formula (\ref{eqn6}) to $f'(s)=g_{n}(\tilde X_s)$,
where $g_n(x)=nx$ for $0\leq x \leq 1/n$ and $g_N(x)=1$ for $x>1/n$.
Since $g_n(\tilde X_u)=0$ for every $u\in ]s, d_s[$ when $s\in{\mathcal J}$, we get
$$\int_{0}^{T_t}g_n(\tilde X_u){\rm d}u\,=\,\int_0^t g_n(X_u){\rm d}u\,.
$$
Letting $n$ tend to $\infty$, we deduce by monotone convergence that
$$\int_{0}^{T_t}{\bf 1}_{\{\tilde X_u>0\}}{\rm d}u\,=\,\int_0^t {\bf 1}_{\{X_u>0\}}{\rm d}u\,.
$$
Since the left-hand side equals $t$, this shows that
(\ref{Cond3}) holds.

Finally, note that for every $t\geq 0$, there is the identity
$$\zeta_t-t :=\inf\{u\geq 0: X_{t+u}=0\}=\inf\{u\geq 0: \tilde X_{T_t+u}=0\}\,.$$
More precisely, either $X_t=0$ and then the two quantities above must be zero,
or $X_t=\tilde X_{T_t}>0$ and then it is immediately seen that
$T_{t+u}=T_{t}+u$ for every $0 \leq u < \zeta_t-t $.
That $(X,V)$ verifies the requirement (\ref{Cond4}) now readily stems  from
the fact that if $S$ is an $({\mathcal F}_t)$-stopping time, then
$T_S$ is a stopping time in the natural filtration of the Wiener process $W$ and
 the strong Markov property of $W$ applied at time $T_S$, as  $X_{t+S}=\tilde X_{T_S+t}$ for every $0\leq t < \zeta_S-S$.
 
 We then observe the identity
 \begin{equation}\label{eqn7}
 V_t=W\circ T_t\,,\qquad \hbox{for all $t\geq 0$, a.s.}
 \end{equation}
 Indeed, in the case when $X_t>0$, we have $V_t=\tilde V_{T_t}=W_{T_t}$.
 In the case when $X_t=0$, since by (\ref{Cond3}), $X$ cannot stay
at $0$ in any neighborhood of $t$, $V$ must take strictly positive values at some times arbitrarily closed to $t$. Because $W$ is continuous and the time-change $T$ right-continuous, this forces $W_{T_{t}}=0$. By (\ref{Cond2}), we conclude that (\ref{eqn7})
also holds in that case.

The regeneration property (\ref{Cond5}) follows from an argument similar to that for (\ref{Cond4}).
Let $S$ be a stopping time in the natural filtration of $X$ with $X_S=0$ a.s.
It is easily seen that $T_{S}$ is a stopping time for the Wiener process $W$,
and by (\ref{Cond2}) and (\ref{eqn7}),  $V_{S}=W_{T_{S}}=0$.   By the strong Markov property,
$W'=(W'_t:=W_{T_{S}+t},t\geq 0)$ is thus a standard Wiener process, which is independent
of the stopped Wiener process $(W_{u\wedge T_{S}}, u\geq 0)$.
Let $Y'$ denote the free Langevin process started from the location $0$ with zero velocity
and driven by $W'$, and by $I'$ its infimum process.
As the free Langevin process $Y$ coincides with its infimum $I$
at time $T_{S}$, we have the identities
$$Y'_t=Y_{T_{S+t}}-Y_{T_{S}}\quad , \quad I'_t=I_{T_{S+t}}-Y_{T_{S}}$$
for every $t\geq 0$. It should now be plain that
$X'=(X'_t:=X_{S+t}, t\geq 0)$ is the reflected Langevin process constructed from $W'$,
which establishes the regeneration property (\ref{Cond5}).
\QED

We now conclude this section by pointing at an interesting connection
between the reflecting Langevin process and a reflecting stable process,
which completes Lemma \ref{Lemma1}. Recall that $L$ denotes the local time at $0$
of the Wiener process and introduce
$$\lambda_t:=L\circ T_t\,,\qquad t\geq 0\,.$$

\begin{proposition}\label{Prop1} {\rm (i)} The process $\lambda=(\lambda_t, t\geq 0)$ serves as a local time at $0$ for the velocity $V$ of the reflecting Langevin process $X$, in the sense that
$$\lambda_t=\lim_{\varepsilon\to 0+}\frac{1}{\varepsilon}\int_0^t
{\bf 1}_{\{0<V_s<\varepsilon\}}{\rm d}s\,,Ê\qquad t\geq 0\,.$$
Moreover the non-decreasing process $\lambda=(\lambda_t, t\geq 0)$ has continuous paths and the support of the Stieltjes measure ${\rm d}\lambda_t$ coincides with
$\{t\geq 0: V_t=0\}$ a.s.

\noindent  {\rm (ii)} Introduce the right-continuous inverse 
$$\lambda^{-1}_t:=\inf\{s\geq 0: \lambda_s>t\}\,,\qquad t\geq 0\,.$$
The time-changed process $X\circ \lambda^{-1}$ can be expressed as
$X\circ \lambda^{-1}_t=\sigma_t-\iota_t$, where
$\sigma$ is the symmetric stable process which appears in Lemma \ref{Lemma1} and
$\iota_t=\inf_{0\leq s\leq t}\sigma_s$. In other words, 
$X\circ \lambda^{-1}$ is a symmetric stable process with index $1/3$ reflected at its infimum.
\end{proposition}

\proof (i) Introduce for every $\varepsilon>0$ and $t\geq 0$
$$L^{(\varepsilon)}_t:=\frac{1}{\varepsilon}\int_0^t{\bf 1}_{\{0<W_s<\varepsilon\}}{\rm d}s\,.$$
The change of variables formula (\ref{eqn6}) yields
$$L^{(\varepsilon)}_{T_t}:=\frac{1}{\varepsilon}\int_0^t{\bf 1}_{\{0<W_{T_s}<\varepsilon\}}{\rm d}s\,.$$
Then recall from (\ref{eqn7}) that $W_{T_s}=V_s$ for all $s\geq 0$, so
taking the limit as $\varepsilon\to 0+$ establishes the first assertion.
Further, as $L^{(\varepsilon)}_{\cdot}$ converges to $L_{\cdot}$ uniformly on every compact time-interval (this follows from the fact that the occupation densities of the Wiener process are jointly continuous), we see that $t\to L_{T_t}$ is continuous a.s. The last assertion can then be proved by a routine argument.

\noindent (ii) Introduce $A_t=\int_0^t{\bf 1}_{\{\tilde X_s>0\}}{\rm d}s$, so that $T_{\cdot}$ is the right-continuous inverse of $A_{\cdot}$. Then $\lambda^{-1}_t=A\circ \tau_t$
and thus $X\circ \lambda^{-1}=\tilde X\circ T\circ A\circ \tau$.
It is easily checked from the strong Markov property of $W$ at $\tau_t$ that $T\circ A\circ \tau_t=\tau_t$ a.s. for each $t\geq0$, 
and as both processes are c\`adl\`ag, the identity holds simultaneously for all $t\geq 0$, a.s. We conclude that
$X\circ \lambda^{-1}_t=Y\circ \tau_t-I\circ {\tau_t}=\sigma_t-\iota_t$. 
\QED
\end{section}

\begin{section}{Uniqueness in distribution}
The purpose of this section is to show that
the requirements (\ref{Cond1}-\ref{Cond5}) 
characterize the law of the reflecting Langevin process.
 Our approach is based on excursion theory  for Markov processes
(see e.g.  Blumenthal \cite{Bl} for details and references);
 let us now explain roughly the main steps.

We shall denote here by $(X,V)$ any c\`adl\`ag process which fulfills
the conditions (\ref{Cond1}-\ref{Cond5}) of Section 2; recall from (\ref{Cond4}) and (\ref{Cond5}) that
the latter has the strong Markov property, and is regenerated at
every stopping time $S$ at which $X_S=0$ a.s.
It is easy to deduce from the right-continuity of the sample paths
and elementary properties of the free Langevin process 
that $X$ returns to $0$ at arbitrarily small times a.s. 
(recall from (\ref{Cond2}) that necessarily $V=0$ at such times); one says that the point $(0,0)$ is regular for $(X,V)$.
This enables us to construct a {\it local time process}
$\ell=(\ell_t, t\geq 0)$, that is a continuous non-decreasing additive functional
which increases exactly on the set of times at which $X$  visits $0$, in the sense that
the support of the random measure ${\rm d} \ell_t$ coincides
with ${\mathcal Z}:=\{t\geq 0: X_t=0\}$ a.s. The local time process
is unique up to some deterministic factor;
the choice of this multiplicative constant is unimportant and just a matter of convention.
{\bf In the sequel, we shall denote generically by $c$ or $c'$ a positive finite constant
whose value may be different in different expressions.}

Next, introduce the (right-continuous) inverse local time
$$\ell^{-1}_t:=\inf\{s\geq 0: \ell_s>t\}\,, \qquad t\geq 0\,.$$
The jumps of $\ell^{-1}$ correspond to the time-intervals on which
$X$ makes an excursion away from $0$, in the sense that
the canonical decomposition of the random open set ${\mathcal Z}^{c}$
on which $X>0$ into interval components, is given by
$${\mathcal Z}^{c}\,=\,\bigcup]\ell^{-1}_{t-},\ell^{-1}_t[\,,$$
where the union is taken over the set of times $t$ at which $\ell^{-1}$ jumps.
The point process 
$$t\to \left(X_{(\ell^{-1}_{t-}+s)\wedge \ell^{-1}_{t}}, s\geq0\right)\,,
\qquad t \hbox{ jump time of }\ell^{-1}\,,$$
is a Poisson point process.
It takes values in the space of excursions, viz. continuous paths in $[0,\infty[$
which
start from $0$, immediately leave $0$, and are absorbed at $0$ after their first return.
Its intensity ${\bf n}$ is called
the excursion measure (of $X$ away from $0$).
The measure ${\bf n}$ does not need to be finite; however 
it always assigns a finite mass to the set of excursion with height greater than 
$\eta$, for any $\eta>0$.
As, by (\ref{Cond3}), $X$ spends no time at $0$, one 
can reconstruct $X$ from the excursion process (this is 
known as It\^o's program), and thus the law of $X$
is determined by the excursion measure ${\bf n}$.
The latter is only defined up to a multiplicative constant
depending on the normalization of the local time.
Thus, in order to establish that the requirements (\ref{Cond1}-\ref{Cond5}) 
characterize the law of $X$, it suffices to show that these
conditions specify ${\bf n}$ up to some deterministic factor, which is the
main goal of this section.

In this direction, we shall denote by $e=(e_t, t\geq 0)$ a generic excursion of
$X$ away from $0$, by $\zeta=\inf\{t>0: e_t=0\}$ the lifetime of the excursion,
and by $v=(v_t, t\geq 0)$ its velocity, so
$$e_t=\int_0^t v_s {\rm d}s\,, \qquad t\geq 0\,.$$
The first  step in the analysis of the excursion measure ${\bf n}$ concerns the behavior of the velocity immediately after the initial time.
More precisely, we are interested in the first instant 
at which the velocity $v$ is $0$ again, viz
$$d:=\inf\{t>0: v_t=0\}\,.$$
The requirements $v_0=0$ and $e_t>0$ for every $0<t<\zeta$
might suggest that the velocity could remain strictly positive immediately
after the initial time, that is that $d>0$. However we shall see that
this intuition is actually wrong.

\begin{lemma}\label{Lemma3} We have
${\bf n}(d>0)=0$.
\end{lemma}
\proof Roughly, the  reason why the velocity cannot remain strictly positive immediately
after the initial time under the excursion measure ${\bf n}$, is that otherwise
the duration $\zeta$ of the excursion $e$ would be too long to allow the reconstruction of the process $X$ from its excursions;
in other words the conditions for It\^o's program would fail. More precisely, recall that the distribution of the duration $\zeta$ of the generic excursion must fulfill
\begin{equation}\label{eqn8}
\int_0^{\infty}(t\wedge 1){\bf n}(\zeta\in {\rm d}t) \, <\, \infty \,;
\end{equation}
 see e.g. Condition (iii) on page 133 in \cite{Bl}.
We shall show that (\ref{eqn8}) can only hold if ${\bf n}(d>0)=0$. 

We suppose hereafter that ${\bf n}(d>0)>0$. It is well-known  that the strong Markov property of $(X,V)$ can be shifted to $(e,v)$ under the excursion measure ${\bf n}$.
It  follows from (\ref{Cond4}) that under the truncated measure ${\bf 1}_{\{d>0\}}{\bf n}$,  the process
$(v_{t\wedge d}, t> 0)$ is Markovian with semigroup given by that of the Wiener process 
in $[0,\infty[$ stopped at its first hitting time of $0$. As  under ${\bf n}$,
the velocity process $v=(v_t, t \geq 0)$ is right-continuous (in fact, continuous on $[0,\zeta[$)  with $v_0=0$, it is well-known that 
this entails that the distribution of $(v_{t\wedge d}, t\geq 0)$ under ${\bf 1}_{\{d>0\}}{\bf n}$
must  be proportional to the It\^o measure   of positive Brownian excursions.

It then follows from the scaling property of It\^o measure  that
for some constant $c>0$,
\begin{equation}\label{levy}
{\bf n}(e_d\in {\rm d}x)=c x^{-4/3}{\rm d}x\,,\qquad x>0\,,
\end{equation}
see Biane and Yor \cite{BY}. Further, as $d$ is a stopping time, the strong Markov property
entails that conditionally on $e_d=x>0$,
 the evolution of $e$ after time $d$  is given
by that of the  Langevin process $Y^{(x)}$ started from $x$ with zero velocity and stopped at $0$, that is at time
$\zeta^{(x)}:=\inf\{t\geq 0: Y^{(x)}_t=0\}$.

It is easily seen from the scaling property of the Wiener process that there is the identity in distribution
$$(Y^{(x)}_t, t\geq 0)\stackrel{\mathcal L}{=}( x Y^{(1)}_{x^{-2/3}t}, t\geq 0)\,,$$
and as a consequence
$$
\zeta^{(x)}\stackrel{\mathcal L}{=} x^{2/3}\zeta^{(1)}\,.
$$
Using (\ref{levy}), we now see that
\begin{eqnarray*}
\int_0^{\infty}(t\wedge 1){\bf n}(\zeta\in {\rm d}t) &\geq& c\int_0^{\infty}x^{-4/3}\E(\zeta^{(x)}\wedge 1) {\rm d}x \\
&=& c\int_0^{\infty}x^{-2/3}\E(\zeta^{(1)}\wedge x^{-2/3}) {\rm d}x\,.
\end{eqnarray*}
The exact distribution of $\zeta^{(1)}$ has been determined by
Goldman \cite{Goldman}, see also Lachal \cite{La91}. We recall 
that the asymptotic behavior of its tail distribution is given by
$$
\P(\zeta^{(1)}>t)\sim c' t^{-1/4}\,, \qquad t\to \infty\,,
$$
where $c'>0$ is some constant, see Proposition 2 in \cite{Goldman}. As a consequence, we have
$$\E(\zeta^{(1)}\wedge x^{-2/3})\geq \int_0^{x^{-2/3}}\P(\zeta^{(1)}>t){\rm d}t
\sim c'   \int_0^{x^{-2/3}}t^{-1/4}{\rm d}t = \frac{4}{3}c'x^{-1/2}\,,\quad x\to 0+\,.
$$
This shows that (\ref{eqn8}) fails and thus ${\bf n}(d>0)$ must be zero.
\QED

The next step of our analysis is provided by Rogozin's solution of the two-sided exit problem for stable L\'evy processes, that we now specify for  the symmetric stable process with index $1/3$ $\sigma$ which arises in Lemma \ref{Lemma1}.
Recall that single points are polar for $\sigma$, so the latter always exits from an interval by a jump.

\begin{lemma}\label{Lemma4}
{\rm (Rogozin \cite{Rogozin})} For every $\varepsilon>0$ and  $x\in]0,\varepsilon[$,
write
$$\varrho_{x,\varepsilon}:=\inf\{t\geq 0: x+\sigma_t\not\in[0,\varepsilon]\}\,.$$
Then   the following identity between sub-probability measures on $[\varepsilon,\infty[$, holds :
$$\P(x+\sigma_{\varrho_{x,\varepsilon}}\in{\rm d}y)
= \frac{1}{2\pi} x^{1/6}(\varepsilon-x)^{1/6}(y-\varepsilon)^{-1/6} y^{-1/6} (y-x)^{-1}
{\rm d}y\,.$$
In particular 
$$\P\left(x+\sigma\hbox{ exits from }[0,\varepsilon]\hbox{ at the upper boundary}\right)
\geq c^{-1} (x/\varepsilon)^{1/6}\,,$$
where $c>0$ is some constant.

\end{lemma}

The result of Rogozin combined with Lemma \ref{Lemma1} enables us to determine a key distribution under the excursion measure ${\bf n}$. 
Introduce for every $\eta >0$  the  passage time
$$
\rho_{\eta}:=\inf\{t\geq 0: e_t\geq \eta\hbox{ and }v_t=0\}\,,$$
with the usual convention that $\inf\emptyset =0$.
It should be plain that
$$ \rho_{\eta}<\infty \Leftrightarrow
\max_{0\leq t \leq \zeta} e_t\geq \eta\,,$$
and in this case, $\rho_{\eta}$ is simply the first instant 
after the first passage time of $e$ at $\eta$ at which the velocity $v$ vanishes.
Observe that $\eta\to \rho_{\eta}$ is non-decreasing,  and also from Lemma \ref{Lemma3} that
\begin{equation}\label{eqn10}
{\bf n}(\rho_{0+}>0)=0\,,\quad \hbox{where} \quad \rho_{0+}=\lim_{\eta\to 0+}
\rho_{\eta}\,.
\end{equation}

 \begin{corollary}\label{Cor1}
There is a finite constant $c>0$ such that for every $\varepsilon >0$,
$${\bf n}(e_{\rho_{\varepsilon}}\in{\rm d}y, \rho_{\varepsilon}<\infty)=
c \varepsilon^{1/6}(y-\varepsilon)^{-1/6}y^{-7/6}{\rm d}y\,,\qquad y\geq \varepsilon\,.$$
\end{corollary}

\proof For the sake of simplicity, we shall assume that $\varepsilon <1$. Fix $\eta>0$.
By standard excursion theory, the condition (\ref{Cond4}) implies that the distribution of
the shifted process
$(e_{\rho_{\eta}+t}, t\geq 0)$ 
under the conditional law ${\bf n}(\cdot \mid  e_{\rho_{\eta}}=x, \rho_{\eta} <\infty)$
is that of a Langevin process started from the location $x$ with zero velocity and stopped at $0$. It follows from Lemmas  \ref{Lemma1} and  \ref{Lemma4}
that for every $x>0$,
$${\bf n}(\rho_1<\infty\mid  e_{\rho_{\eta}}=x,\rho_{\eta} <\infty)\,\geq \, c^{-1} (1\wedge x)^{1/6}\,,$$
and therefore
$$\int_{]0,\infty[} (1\wedge x)^{1/6}{\bf n}(e_{\rho_{\eta}}\in{\rm d}x, \rho_{\eta} <\infty)\,\leq \, c'\,.$$
Thus 
$$\mu_{\eta}({\rm d}x):=
(1\wedge x)^{1/6}
{\bf n}(e_{\rho_{\eta}}\in{\rm d}x, \rho_{\eta}<\infty)\,,\qquad{\eta>0}$$ is a family of bounded measures on $\R_+$, and we deduce from (\ref{eqn10}) that
the following limit holds in the sense of weak convergence of finite measures on $\R_+$
as $\eta$ tends to $0$ along some sequence, say $(\eta_k, k\in\N)$ :
\begin{equation}\label{eqn11}
\lim_{k\to\infty} \mu_{\eta_k}({\rm d}x)=a\delta_0({\rm d} x),
\end{equation}
where $a\geq 0$ is some constant.

Next, we observe from the strong Markov property (and again Lemma \ref{Lemma1})
that the measure ${\bf n}(e_{\rho_{\epsilon}}\in{\rm d}y, \rho_{\epsilon}<\infty)$ on $]\varepsilon,\infty[$ can be expressed in the form
$${\bf 1}_{\{y\geq \varepsilon\}}{\bf n}(e_{\rho_{\eta}}\in{\rm d}y, \rho_{\eta}<\infty)+
\int_{]0,\varepsilon[}{\bf n}(e_{\rho_{\eta}}\in{\rm d}x, \rho_{\eta}<\infty)
\P(x+\sigma_{\varrho_{x,\varepsilon}}\in{\rm d}y)\,.$$

On the one hand, (\ref{eqn10}) entails that the first term in the sum converges to
$0$ as $\eta\to 0+$. On the other hand, we may rewrite the second term
in the form
\begin{eqnarray*}
& &\int_{]0,\varepsilon[}
\mu_{\eta}({\rm d}x) x^{-1/6}
\P(x+\sigma_{\varrho_{x,\varepsilon}}\in{\rm d}y)\\
&=& \frac{1}{2\pi}(y-\varepsilon)^{-1/6} y^{-1/6}\left(
\int_{]0,\varepsilon[}
\mu_{\eta}({\rm d}x) (\varepsilon-x)^{1/6} (y-x)^{-1}\right)
{\rm d}y\,,
\end{eqnarray*}
where the equality stems from
 Lemma \ref{Lemma4}.
Taking the limit as $\eta\to 0+$, we conclude from (\ref{eqn11}) that
$${\bf n}(e_{\rho_{\epsilon}}\in{\rm d}y, \rho_{\epsilon}<\infty)
=\frac{a}{2\pi} \varepsilon^{1/6}(y-\varepsilon)^{-1/6} y^{-7/6}{\rm d}y\,,\qquad y>\varepsilon\,,$$
which is our  claim.  \QED

Roughly speaking, Corollary \ref{Cor1} provides a substitute for the entrance law under the excursion measure.
We are now able to establish the uniqueness in distribution for the reflecting
Langevin process.

\begin{theorem}\label{Thm2} Any c\`adl\`ag process which fulfills the conditions
(\ref{Cond1}-\ref{Cond5}) is distributed as the process constructed in Section 2.
\end{theorem}

\proof Recall that
all that is needed is to check that the conditions
(\ref{Cond1}-\ref{Cond5}) characterize the excursion measure ${\bf n}$ up to some deterministic factor.

In this direction, it is convenient to introduce
the space ${\mathcal C}_b$ of bounded continuous paths $\omega: \R_+\to \R$
 endowed with the supremum norm, and to denote by
$\theta_t: {\mathcal C}_b\to {\mathcal C}_b$ the usual shift operator.  We write
${\P}_x^{\partial}$ for the probability measure on ${\mathcal C}_b$ induced by the Langevin process started from $x>0$
with zero velocity and stopped at $0$.
 For every $\eta>0$ and every path $\omega$, set 
$$\varrho_{\eta}:=\inf\{t\geq 0: \bar\omega(t)>\eta, \hbox{ and }\bar\omega(t)>\omega(t)\}\,,$$
  where $\bar\omega(t):=\max_{0\leq s \leq t}\omega(s)$.
 Observe that for ${\bf n}$-almost every paths, $\rho_{\eta}=\varrho_{\eta}$.
 Recall also (\ref{eqn10})
 
Consider a  continuous bounded
functional $F: {\mathcal C}_b\to \R$  which is identically $0$
 on some neighborhood of the degenerate path  $\omega\equiv 0$.
 Then for ${\bf n}$-almost every paths $\omega$, we have
 $$\lim_{\eta\to 0+}{\bf 1}_{\{\varrho_{\eta}<\infty\}}F(\theta_{\varrho_{\eta}}(\omega))=F(\omega)\,,$$
  and by dominated convergence
 $${\bf n}(F(\omega))=\lim_{\eta\to 0+}{\bf n}(F(\theta_{\varrho_{\eta}}(\omega)), \varrho_{\eta}<\infty)\,.$$
 Since, by the strong Markov property and Corollary \ref{Cor1},
 $${\bf n}(F(\theta_{\varrho_{\eta}}(\omega)),\varrho_{\eta}<\infty)=c  \eta^{1/6}
 \int_{\eta}^{\infty} \E_x^{\partial}(F(\omega))
(x-\eta)^{-1/6}x^{-7/6}{\rm d}x\,,$$
this determines ${\bf n}$ up to a multiplicative constant. \QED
 
 We point out that a slight variation of the argument in the proof of Theorem \ref{Thm2} 
 shows that  for every continuous bounded
functional $F: {\mathcal C}_b\to \R$  which is identically $0$
 on some neighborhood of $\omega\equiv 0$, there is the approximation
\begin{equation}\label{eqn12}
{\bf n}(F(e))\,=\, \lim_{x\to 0+} cx^{-1/6}\E_x^{\partial}(F(\omega))\,.
\end{equation}

\end{section}

\begin{section}{Scaling exponents}
We now conclude this paper by answering some natural questions about the scaling exponents of the reflecting Langevin process.

\begin{proposition}\label{Prop2} {\rm (i) } For every $a>0$, there is the identity in distribution
$$(a^{-3/2}X_{at}, a^{-1/2}V_{at})_{t\geq 0} \stackrel{\mathcal L}{=}(X_t,V_t)_{t\geq 0}\,.
$$

\noindent {\rm (ii) } For every $a>0$, the distribution of
$(a^{-3/2}e_{at})_{t\geq 0}$ under  the excursion measure ${\bf n}$
is $a^{-1/4}{\bf n}$.

\end{proposition}

\proof (i) It is immediately seen from the scaling property of the Wiener process
that the free Langevin process $Y$ started from $0$ with zero initial velocity
fulfills
$$(a^{-3/2}Y_{at}, a^{-1/2}\dot Y_{at})_{t\geq 0} \stackrel{\mathcal L}{=}(Y_t,\dot Y_t)_{t\geq 0}\,.$$
The first assertion now follows readily, as the time substitution $(T_t)_{t\geq 0}$ does not affect the scaling property.

\noindent (ii) Standard arguments of excursion theory combined with the scaling property (i) show that the distribution of
$(a^{-3/2}e_{at})_{t\geq 0}$ under  the excursion measure ${\bf n}$
must be proportional to ${\bf n}$. The value of the factor can be determined using Corollary
\ref{Cor1}. Indeed, if we denote by $h:=\max_{t\geq 0}e_t$ the height of the excursion, then we see that 
\begin{equation}\label{eqn13}
{\bf n}(h>x)=cx^{-1/6}\,,Ê\qquad x>0\,.
\end{equation}
 Thus  ${\bf n}(a^{-3/2}h>x)=a^{-1/4}
{\bf n}(h>x)$, which entails our claim. \QED

In the proof of the scaling property of the excursion measure, we have determined the law of the height of the excursion; cf. (\ref{eqn13}). More generally, the scaling property
enables us to specify the law under the excursion measure of any variable
which enjoys the self-similarity property. Here is an example of application.

\begin{corollary}\label{Cor2}  {\rm (i) }
The tail-distribution of the lifetime $\zeta$ of an excursion is given by
$${\bf n}(\zeta>t)=c t^{-1/4}\,,\qquad t>0\,.$$
As a consequence, the inverse local time process $\ell^{-1}$ is a stable subordinator
with index $1/4$, and in particular
the exact  Hausdorff function of the zero set
${\mathcal Z}=\{
t\geq 0: X_t=0\}$ 
of the reflecting Langevin process of is $H(\varepsilon)=\varepsilon^{1/4}(\ln\ln 1/\varepsilon)^{3/4}$  a.s.

\noindent {\rm (ii) } The tail-distribution of the velocity of the excursion immediately before
its lifetime is given by
$${\bf n}(v_{\zeta-}<-x)=c' x^{-1/2}\,,\qquad x>0\,.$$
\end{corollary}

\proof The formulas for the tail-distributions follow immediately from
Proposition \ref{Prop2}. In particular, the first one entails that the inverse local time
$\ell^{-1}$ of the reflecting Langevin process is a stable subordinator with index 
$1/4$. The assertion about the exact Hausdorff measure then
derives from a well-known result due to Taylor and Wendel \cite{TW}. \QED

It may be interesting to compare Corollary \ref{Cor2}(i) with the case when the boundary is reflecting, in the sense that a particle arriving at $0$ with incoming velocity $-v<0$ bounces back with  velocity $v$. It is straightforward to check that if $Y$ is a free Langevin process, then $|Y|$ describes the process that bounces at the boundary. 
The set of times $\{t\geq 0: |Y_t|=0\}$ at which the particle hits the obstacle
is then countable a.s., and its only accumulation point is $t=0$.
We refer to Lachal \cite{La97} for much more on this topic.

To conclude, we also point out that
Corollary \ref{Cor2}(ii) entails that the total energy absorbed by the obstacle at time $t$,
$${\mathcal E}_t:=\frac{1}{2}\sum_{s\in{\mathcal Z}\cap [0,t]}V_{s-}^2\,,$$
is finite a.s. More precisely, the time-changed process ${\mathcal E}\circ \ell^{-1}$
is a stable subordinator with index $1/4$.
\end{section}

\vskip 1cm
\noindent
{\bf Acknowledgment :} This paper has been motivated by a lecture and a series of questions by Bertrand Maury. In particular, Bertrand raised the key problem  of the existence of a reflecting Langevin process
with an absorbing boundary, which he observed in numerical simulations.

\end{document}